\documentclass{article}    
\usepackage{amstext}
\def\qed{$\rlap{$\sqcap$}\sqcup$}
\usepackage{latexsym}

\begin{document}           

\begin{center}
{\ }\\
{\huge {\bf Improving the bounds of the Multiplicity Conjecture: the codimension 3 level case}} \\ [.250in]
{\large FABRIZIO ZANELLO\\
Department of Mathematics, Royal Institute of Technology, 100 44 Stockholm, Sweden
\\E-mail: zanello@math.kth.se}
\end{center}

{\ }\\
\\
ABSTRACT. The Multiplicity Conjecture (MC) of Huneke and Srinivasan provides upper and lower bounds for the multiplicity of a Cohen-Macaulay algebra $A$ in terms of the shifts appearing in the modules of the minimal free resolution (MFR) of $A$. All the examples studied so far have lead to conjecture (see $[HZ]$ and $[MNR2]$) that, moreover, the bounds of the MC are sharp if and only if $A$ has a pure MFR. Therefore, it seems a reasonable - and useful - idea to seek better, if possibly {\it ad hoc}, bounds for particular classes of Cohen-Macaulay algebras.\\
In this work we will only consider the codimension 3 case. In the first part we will stick to the bounds of the MC, and show that they hold for those algebras whose $h$-vector is that of a compressed algebra.\\
In the second part, we will (mainly) focus on the level case: we will construct new conjectural upper and lower bounds for the multiplicity of a codimension 3 level algebra $A$, which can be expressed exclusively in terms of the $h$-vector of $A$, and which are better than (or equal to) those provided by the MC. Also, our bounds can be sharp even when the MFR of $A$ is not pure.\\
Even though proving our bounds still appears too difficult a task in general, we are already able to show them for some interesting classes of codimension 3 level algebras $A$: namely, when $A$ is compressed, or when its $h$-vector $h(A)$ ends with $(...,3,2)$. Also, we will prove our lower bound when $h(A)$ begins with $(1,3,h_2,...)$, where $h_2\leq 4$, and our upper bound when $h(A)$ ends with $(...,h_{c-1},h_c)$, where $h_{c-1}\leq h_c+1$.\\
\\
\\
\section{Introduction}
\indent
{\large

In their 1985 article $[HM]$, Huneke and Miller showed that, when  a Cohen-Macaulay algebra $A$ of codimension $r$ has a pure Minimal Free Resolution (MFR), then the multiplicity of $A$ equals the product of the different shifts appearing in its MFR divided by the factorial of $r$. The attempt to generalize this seminal result, by relating the multiplicity of any Cohen-Macaulay algebra to the shifts appearing in the modules of its MFR, led to the formulation of the so-called {\it Multiplicity Conjecture} (MC), due to Huneke and Srinivasan: namely, for any Cohen-Macaulay algebra $A$, its multiplicity times the factorial of its codimension is bounded from below (respectively, from above) by the product of the smallest (respectively, largest) shifts of the modules of the MFR of $A$.\\\indent
Notice that, since the multiplicity and the MFR of a Cohen-Macaulay algebra $A$ are preserved when one considers the artinian reductions of $A$, it suffices to study the Multiplicity Conjecture for artinian algebras. We just recall here that the MC has been extended to the non-Cohen-Macaulay case by Herzog and Srinivasan, who conjectured (the stronger fact) that the upper bound of the MC actually holds for {\it any} graded algebra (whereas examples have shown that the lower bound in general does not hold when we drop the Cohen-Macaulay hypothesis). However, we will only consider the Cohen-Macaulay case in this article\\\indent
The MC has been attacked by a number of researchers over the last years, but has so far been settled only in particular cases: among them, for codimension 2 algebras (see $[HS]$; see also $[MNR1]$, where the bounds of the MC are improved), Gorenstein algebras of codimension 3 (see $[MNR1]$, where again better bounds are shown to hold), algebras with a quasi-pure resolution ($[HS]$), standard determinantal ideals ($[Mi]$), componentwise linear ideals ($[Ro]$), complete intersections ($[HS]$) and powers thereof ($[GV]$). We refer the reader to the recent works $[Fr]$ and $[FS]$ for a comprehensive history of all the main results obtained to date on the MC.\\\indent
All the results achieved so far have lead to the further conjecture - due to Herzog and Zheng ($[HZ]$) and Migliore, Nagel and R\"omer ($[MNR2]$) - that the bounds provided by the MC are sharp if and only if the algebra we are considering has a pure MFR. Therefore, as Migliore {\it et al.} did successfully in $[MNR1]$ for the two cases mentioned above, it seems suitable to attack the MC by seeking better - if {\it ad hoc} - bounds any time we consider particular classes of algebras (without a pure MFR). We will take this approach here in studying codimension 3 level algebras, for which we will supply better conjectural bounds in the third section of this paper.\\
\\\indent
Throughout this work we will only consider codimension 3 algebras. In the next section, we will prove the MC (without explicit improvements) for those algebras having the same $h$-vector as that of a compressed algebra.\\\indent
In the third section, we will restrict our attention to level algebras $A$, and construct {\it ad hoc} upper and lower bounds, which are better than (or, at worse, equal to) those of the MC; also, they are in general sharp, even when the MFR of $A$ is not pure. These bounds have the possibly great advantage of being entirely recovered from the $h$-vector of $A$, and therefore one does not need explicit information on the MFR's of codimension 3 level algebras, but just on their $h$-vectors, to determine whether the bounds hold. However, since in codimension 3 even the structure of level $h$-vectors is still far from being completely understood, at this point we are able to prove our conjecture only for some interesting special classes of codimension 3 level algebras $A$: namely, when $A$ is compressed, or when its $h$-vector $h(A)$ ends with $(...,3,2)$. Furthermore, we will prove our lower bound when $h(A)$ begins with $(1,3,h_2,...)$, where $h_2\leq 4$, and our upper bound when $h(A)$ ends with $(...,h_{c-1},h_c)$, where $h_{c-1}\leq h_c+1$.\\\indent
Finally, we will show that our bounds cannot in general be extended to the non-level case. However, even though our lower bound does not hold for all codimension 3 algebras (we will supply some non-level counterexamples), we will see that when it is verified for a given algebra $A$, then the lower bound of the MC holds for $A$ as well. We will exploit this fact to easily show the lower bound of the MC for any algebra whose $h$-vector begins with $(1,3,4,5,...)$.\\
\\\indent
Let us now fix the main definitions we will need in this paper. As we said, since we are only studying the MC in the Cohen-Macaulay case, we can suppose without loss of generality that our (standard graded) algebras $A$ are artinian. We set $A=R/I$, where $R=k[x_1,...,x_r]$, $k$ is a field of characteristic zero, $I$ is a homogeneous ideal of $R$ and the $x_i$'s all have degree 1.\\\indent
The {\it $h$-vector} of $A$ is $h(A)=h=(h_0,h_1,...,h_c)$, where $h_i=\dim_k A_i$ and $c$ is the last index such that $\dim_k A_c>0$. Since we may suppose, without loss of generality, that $I$ does not contain non-zero forms of degree 1, $r=h_1$ is defined as the {\it codimension} of $A$. The {\it multiplicity} of $A$ is $e=1+h_1+...+h_c$, that is the dimension of $A$ as a $k$-vector space. The {\it initial degree} of $I$ is the least degree $t$ where $I$ is non-zero, or equivalently, the least index such that $h_t$ is not {\it full-dimensional} (i.e. $h_t<{r-1+t\choose t}$).\\\indent 
The {\it socle} of $A$ is the annihilator of the maximal homogeneous ideal $\overline{m}=(\overline{x_1},...,\overline{x_r})\subseteq A$, namely soc$(A)=\lbrace a\in A {\ } \mid {\ } a\overline{m}=0\rbrace $. Since soc$(A)$ is a homogeneous ideal, we define the {\it socle-vector} of $A$ as $s(A)=s=(s_0,s_1,...,s_c)$, where $s_i=\dim_k$soc$(A)_i$. Notice that $h_0=1$, $s_0=0$ and $s_c=h_c>0$. The integer $c$ is called the {\it socle degree} of $A$ (or of $h$). The {\it type} of the socle-vector $s$ (or of the algebra $A$) is type$(s)=\sum_{i=0}^cs_i$.\\\indent 
If $s=(0,0,...,0,s_c)$, we say that the algebra $A$ is {\it level} (of type $s_c$). In particular, if $s_c=1$, $A$ is {\it Gorenstein}. With a slight abuse of notation, we will sometimes refer to an $h$-vector as Gorenstein (or level) if it is the $h$-vector of a Gorenstein (or level) algebra.\\\indent
The {\it minimal free resolution} ({\it MFR}) of an artinian algebra $A$ is an exact sequence of $R$-modules of the form: $$0\longrightarrow F_r\longrightarrow F_{r-1}\longrightarrow ...\longrightarrow F_1\longrightarrow R\longrightarrow R/I \longrightarrow 0,$$
where, for $i=1,...,r$, $$F_i=\bigoplus_{j=m_i}^{M_i}R^{\beta_{i,j}}(-j),$$  and all the homomorphisms have degree 0.\\\indent The $\beta_{i,j}$'s are called the {\it graded Betti numbers} of $A$.\\\indent 
Then $\beta_{1,j}$ is the number of generators of degree $j$ of $I$. It is well-known that $F_r=\oplus_{j=1}^cR^{s_j}(-j-r)\neq 0$. Hence, the socle-vector may also be computed by considering the graded Betti numbers of the last module of the MFR. In particular, an artinian algebra $A$ is level of socle degree $c$ and type $s_c$ if and only if $F_r=R^{s_c}(-c-r)$. The MFR of an algebra $A$ is called {\it pure} if it has only one different shift in each module.\\\indent 
Let $h(z)=\sum_{i=0}^ch_iz^i$ be the {\it Hilbert series} of $A$ (note that, since $A$ is artinian, $h(z)$ here is in fact a polynomial). In particular, we have $e=h(1)$. The MFR and the Hilbert series of $A$ are related by the following well-known formula (e.g., see $[FL]$, p. 131, point (j) for a proof): \begin{equation}\label{sh}h(z)(1-z)^r=1+\sum_{i,j}(-1)^i\beta_{i,j}z^j.\end{equation}\indent
We are now ready to state the Huneke-Srinivasan Multiplicity Conjecture:\\
\\\indent
{\bf Conjecture 1.1 (Multiplicity Conjecture).} {\it Let $A$ be an artinian algebra of codimension $r$ and multiplicity $e$, and let, for $i=1,2,...,r$, $m_i$ and $M_i$ be, respectively, the smallest and the largest shift appearing in the $i$-th module $F_i$ of the MFR of $A$. Then:
$${m_1\cdot m_2\cdot \cdot \cdot m_r\over r!}\leq e\leq {M_1\cdot M_2\cdot \cdot \cdot M_r\over r!}.$$}
\\
\section{Compressed $h$-vectors of codimension 3}
\indent

The purpose of this section is to show that all the codimension 3 algebras whose $h$-vector is the $h$-vector of a compressed algebra satisfy the Multiplicity Conjecture. The idea of a compressed algebra is a natural concept which first appeared (for the Gorenstein case) in Emsalem-Iarrobino's 1978 seminal paper $[EI]$, and describes those algebras having the (entry by entry) maximal $h$-vector among all the algebras with given codimension and socle-vector. Compressed algebras and their $h$-vectors were extensively studied in the Eighties by Iarrobino ($[Ia]$) and Fr\"oberg-Laksov ($[FL]$) - who restricted their attention to the very natural case where the socle-vectors have \lq \lq enough" initial entries equal to 0; see below for the exact definition -, and recently by this author in full generality - see $[Za1]$ and $[Za2]$, where we have defined {\it generalized} compressed algebras.\\\indent
Let us now fix a codimension $r$ and a socle-vector $s=(s_0=0,s_1,...,s_c)$.\\
\\\indent
{\bf Definition-Remark 2.1.} Following $[FL]$, define, for $d=0,1,...,c$, the integers $$r_d=N(r,d)-N(r,0)s_d-N(r,1)s_{d+1}-...-N(r,c-d)s_c,$$ where we set $$N(r,d)=\dim_kR_d={r-1+d \choose d}.$$\indent
It is easy to show (cf. $[FL]$) that $r_0<0$, $r_c\geq 0$ and $r_{d+1}>r_d$ for every $d$. Define $b$, then, as the unique index such that $1\leq b\leq c$, $r_b\geq 0$ and $r_{b-1}<0$.\\
\\\indent
Let $S=k[y_1,y_2,...,y_r]$, and consider $S$ as a graded $R$-module where the action of $x_i$ on $S$ is partial differentiation with respect to $y_i$. Recall that, in the theory of {\it inverse systems} (for which we refer the reader to $[Ge]$ and $[IK]$), the $R$-submodule $M=I^{-1}$ of $S$, (bijectively) associated to the algebra $R/I$ having socle-vector $s$, is generated by $s_i$ elements of degree $i$, for $i=1,2,...,c$. Furthermore, the $h$-vector of $R/I$ is given by the number of linearly independent partial derivatives obtained in each degree by differentiating the generators of $M$.\\\indent
The number $$N(r,d)-r_d=N(r,0)s_d+N(r,1)s_{d+1}+...+N(r,c-d)s_c$$ is therefore an upper bound for the number of linearly independent derivatives supplied in degree $d$ by the generators of $M$ and, hence, is also an upper bound for the $h$-vector of $R/I$. This is the reason for the introduction of the numbers $r_d$.\\
\\\indent
{\bf Proposition 2.2} (Fr\"oberg-Laksov). {\it Fix a codimension $r$ and a socle-vector $s=(0,s_1,...,s_c)$. Then an upper bound for the $h$-vectors of all the algebras having data $(r,s)$ is given by $$H=(h_0,h_1,...,h_c),$$ where, for $i=0,1,...,c$, $$h_i=\min \lbrace N(r,i)-r_i,N(r,i)\rbrace .$$}
\\\indent
{\bf Proof.} See $[FL]$, Proposition 4, i). (Fr\"oberg and Laksov gave a direct proof of this proposition; notice that a second proof immediately follows from our comment on inverse systems and the numbers $r_d$. The same upper bound was already supplied by Iarrobino ($[Ia]$) under the natural restriction $s_1=...=s_{b-1}=0$.){\ }{\ }\qed \\
\\\indent
{\bf Theorem 2.3} (Iarrobino, Fr\"oberg-Laksov). {\it Let $r$ and $s$ be as above. If $s_1=...=s_{b-1}=0$, then the upper bound $H$ of Proposition 2.2 is actually the $h$-vector of \lq \lq almost all" the algebras (that is, those parameterized by a suitable non-empty Zariski-open set) having data $(r,s)$.}\\
\\\indent
{\bf Proof.} See $[Ia]$, Theorem II A; $[FL]$, Proposition 4, iv) and Theorem 14.{\ }{\ }\qed \\
\\\indent
{\bf Definition 2.4} (Iarrobino). Fix a pair $(r,s)$ such that $s_1=...=s_{b-1}=0$. An algebra having data $(r,s)$ is called {\it compressed (with respect to the pair $(r,s)$)} if its $h$-vector is the upper bound $H$ of Proposition 2.2.\\
\\\indent
Thus, Theorem 2.3 shows the existence of compressed algebras and provides an explicit description of their $h$-vectors. We are now going to see that the MFR's of compressed algebras also have a very nice shape; that is, in each module (of course, except possibly for the last one, which represents the socle) they have at most two different shifts. Precisely:\\
\\\indent
{\bf Proposition 2.5} (Fr\"oberg-Laksov). {\it Let $A$ be a compressed algebra (with respect to a given pair $(r,s)$ such that $s_1=...=s_{b-1}=0$). Then, for each $i=1,2,...,r-1$, the $i$-th module $F_i$ of the MFR of $A$ has at most two different shifts, occurring in degrees $t+i-1$ and $t+i$.}\\
\\\indent
{\bf Proof.} See $[FL]$, Proposition 16.{\ }{\ }\qed \\
\\\indent
Let us now restrict our attention to codimension $r=3$. Next we show that the Multiplicity Conjecture holds for all compressed algebras; we will see later that this result implies the MC for any algebra whose $h$-vector is the same as that of a compressed algebra.\\
\\\indent
{\bf Theorem 2.6.} {\it Let $A$ be a codimension 3 compressed algebra. Then the MC holds for $A$.}\\
\\\indent
{\bf Proof.} With the notation above, let $h$ be the $h$-vector of a compressed algebra $A=R/I$ of codimension 3 and socle-vector $s=(0,0,...,0,s_q,...,s_c)$, where we suppose that $q$ is the smallest index such that $s_q>0$. Hence, by definition, $q\geq b$. Notice that the initial degree $t$ of $I$ equals $b$ if $r_b>0$, whereas $t=b+1$ if $r_b=0$. Also, by definition, we have \begin{equation}\label{t}N(3,c-t)s_c+N(3,c-t-1)s_{c-1}+...+N(3,1)s_{t+1}+N(3,0)s_t<N(3,t).\end{equation}\indent
By Proposition 2.5, the MFR of $A$ has the form:
$$ 0\longrightarrow F_3=\oplus_{j=q}^cR^{s_j}(-(j+3))\longrightarrow F_{2}=R^{\beta_{2,t+1}}(-(t+1))\oplus R^{\beta_{2,t+2}}(-(t+2))$$$$\longrightarrow F_1=R^{\beta_{1,t}}(-t)\oplus R^{\beta_{1,t+1}}(-(t+1))
\longrightarrow R\longrightarrow R/I \longrightarrow 0.$$\indent
The multiplicity of $A$ is
$$e=1+N(3,1)+...+N(3,t-1)+(1+N(3,1)+...+N(3,c-t))s_c$$$$+(1+N(3,1)+...+N(3,c-t-1))s_{c-1}+...+(1+N(3,1)+...+N(3,q-t))s_{q}$$\begin{equation}\label{e}=N(4,t-1)+N(4,c-t)s_c+N(4,c-t-1)s_{c-1}+...+N(4,q-t)s_{q},\end{equation}
the last equality following from the combinatorial identity $\sum_{i=0}^dN(r,i)=N(r+1,d)$.\\\indent
In order to see which are the smaller and the larger shifts in $F_1$ and $F_2$, we need to determine when $\beta_{2,t+2}$ is positive and to distinguish the three cases $\beta_{2,t+1}-\beta_{1,t+1}< 0$, $\beta_{2,t+1}-\beta_{1,t+1}> 0$ and $\beta_{2,t+1}-\beta_{1,t+1}= 0$.\\\indent
By formula (\ref{sh}), since $\beta_{3,t+2}=s_{t-1}$, we easily have
$$\beta_{2,t+2}-s_{t-1}=h_{t+2}-3h_{t+1}+3h_t-h_{t-1}=$$$$N(3,c-t-2)s_c+N(3,c-t-3)s_{c-1}+...+N(3,0)s_{t+2}-3(N(3,c-t-1)s_c+...+N(3,1)s_{t+2}$$$$+N(3,0)s_{t+1})+3(N(3,c-t)s_c+...+N(3,2)s_{t+2}+N(3,1)s_{t+1}+N(3,0)s_t)-N(3,t-1).$$\indent
From the identity $3N(3,i)-3N(3,i-1)+N(3,i-2)=N(3,i+1)$, it follows that
$$\beta_{2,t+2}-s_{t-1}=N(3,c-t-1)s_c+N(3,c-t)s_{c-1}+...+N(3,2)s_{t+1}+N(3,1)s_t-N(3,t-1)$$$$=-r_{t-1}-s_{t-1}.$$\indent
Thus, $\beta_{2,t+2}=-r_{t-1}$. Therefore, it immediately follows that $\beta_{2,t+2}= 0$ if and only if $b=t-1$, if and only if $r_b=0$ - that is, when $A$ is called {\it extremely} compressed.\\\indent
Again by formula (\ref{sh}), we have:
$$\beta_{2,t+1}-\beta_{1,t+1}=h_{t+1}-3h_t+3h_{t-1}-h_{t-2}=$$$$N(3,c-t-1)s_c+N(3,c-t-2)s_{c-1}+...+N(3,0)s_{t+1}-3(N(3,c-t)s_c+$$$$N(3,c-1-t)s_{c-1}+...+N(3,1)s_{t+1}+N(3,0)s_t)+3N(3,t-1)-N(3,t-2).$$\indent
Hence, from the identity $3N(3,i)-N(3,i-1)=(i+1)(i+3)$, we immediately get
\begin{equation}\label{diff}\beta_{2,t+1}-\beta_{1,t+1}=t(t+2)-\sum_{j=q}^cs_j(j-t+1)(j-t+3).\end{equation}\indent
Let us first consider the case $\beta_{2,t+2}=0$. We have seen that this is equivalent to $r_b=0$. In order to show the MC, we have to prove that \begin{equation}\label{0}{t(t+1)(q+3)\over 6}\leq e\leq {t(t+1)(c+3)\over 6}.\end{equation}\indent
The first inequality is $${t(t+1)(q+3)\over 6}\leq N(4,t-1)+N(4,c-t)s_c+N(4,c-t-1)s_{c-1}+...+N(4,q-t)s_{q}.$$\indent
If we bring the summand $N(4,t-1)$ to the l.h.s., it is easy to see that the previous inequality is a consequence of the following: \begin{equation}\label{ww}t(t+1)\leq (c-t+2)(c-t+3)s_c+(c-t+1)(c-t+2)s_{c-1}+...+(q-t+2)(q-t+3)s_q.\end{equation}\indent
But (\ref{ww}) means $$N(t-1,3)\leq N(3,c-t+1)s_c+N(3,c-t)s_{c-1}+...+N(3,q-t+1)s_q,$$ which is true (actually, an equality holds), since $t-1=b$ and $r_b=0$. This proves the first inequality of (\ref{0}).\\\indent
In order to prove the second inequality, recall that $r_b=0$ means \begin{equation}\label{rb}N(3,t-1)=N(3,c-t+1)s_c+N(3,c-t)s_{c-1}+...+N(3,1)s_t+N(3,0)s_{t-1}.\end{equation}\indent
We want to show that $$N(4,t-1)+N(4,c-t)s_c+N(4,c-t-1)s_{c-1}+...+N(4,0)s_{t}\leq {t(t+1)(c+3)\over 6}.$$\indent
By (\ref{e}), (\ref{rb}), and the identities $N(r,i)-N(r-1,i)=N(r,i-1)$ and $\sum_{i=0}^dN(r,i)=N(r+1,d)$, it is easy to see that the last inequality is equivalent to $$N(4,t-2)+N(4,c-t+1)s_c+N(4,c-t)s_{c-1}+...+N(4,1)s_{t}\leq {t(t+1)(c+3)\over 6}.$$\indent
If we bring $N(4,t-2)$ to the r.h.s., one moment's thought shows that it suffices to prove that
$$N(3,c-t+1)s_c+N(3,c-t)s_{c-1}+...+N(3,1)s_t\leq N(3,t-1).$$\indent
But, since $t-1=b$, this is equivalent to $$-r_b-s_b+N(3,t-1)\leq N(3,t-1),$$ i.e. $-r_b-s_b\leq 0$, which is true since $r_b=0$ and $s_b\geq 0$. This completes the proof of (\ref{0}).\\\indent
Hence, from now on, suppose that $r_b>0$, i.e. that $\beta_{2,t+2}>0$. Thus, $b=t$ and $s_{t-1}=0$. Let us first consider the case where the r.h.s. of formula (\ref{diff}) is lower than 0. Therefore, $\beta_{2,t+1}$ can possibly be equal to 0, whereas we must have $\beta_{1,t+1}>0$. Hence, in order to prove the MC, it suffices to show that
\begin{equation}\label{ii}{t(t+2)(q+3)\over 6}\leq e\leq {(t+1)(t+2)(c+3)\over 6}.\end{equation}\indent
We will actually show (\ref{ii}) supposing that the r.h.s. of (\ref{diff}) is less than {\it or equal to} 0. As far as the first inequality of (\ref{ii}) is concerned, what we want to prove is that
$$t(t+1)(t+2)+(c-t+1)(c-t+2)(c-t+3)s_c+(c-t)(c-t+1)(c-t+2)s_{c-1}$$$$+...+(q-t+1)(q-t+2)(q-t+3)s_q\geq t(t+2)(q+3).$$\indent
To this purpose, it suffices to show that $$(c-t+1)(c-t+3)s_c+...+(q-t+1)(q-t+3)s_q\geq t(t+2),$$ which is true under the current assumption that the r.h.s. of formula (\ref{diff}) is less than or equal to 0. This proves the first inequality of (\ref{ii}).\\\indent
As for the second, we want to show that
$$t(t+1)(t+2)+(c-t+1)(c-t+2)(c-t+3)s_c+(c-t)(c-t+1)(c-t+2)s_{c-1}$$$$+...+(q-t+1)(q-t+2)(q-t+3)s_q\leq t(t+2)(c+3).$$\indent
Similarly to above, {\it a fortiori} it is enough to prove that $$(c-t+1)(c-t+2)s_c+...+(q-t+1)(q-t+2)s_q\leq (t+1)(t+2),$$ i.e. that
$$N(3,c-t)s_c+...+N(3,q-t)s_q\leq N(3,t).$$\indent
But the last inequality means exactly $r_t\geq 0$, and this is true since $t=b$. This completes the proof of (\ref{ii}).\\\indent
We now want to show the MC when the r.h.s. of formula (\ref{diff}) is greater than 0 (recall that we are always under the hypothesis $\beta_{2,t+2}>0$). In this case, $\beta_{1,t+1}$ can happen to be 0, whereas we always have $\beta_{2,t+1}>0$. Thus, we want to show that
\begin{equation}\label{iii}{t(t+1)(q+3)\over 6}\leq e\leq {t(t+2)(c+3)\over 6}.\end{equation}\indent
We omit the computations here, since they are just \lq \lq symmetric" to those performed for the previous case, and again also hold when the r.h.s. of formula (\ref{diff}) is equal to 0. We just remark that the inequality on the r.h.s. of (\ref{diff}) is employed only in proving the second inequality of (\ref{iii}), and that, at the end, one uses the fact that $t=b$ implies $r_{t-1}< 0$ to show the first inequality of (\ref{iii}).\\\indent
Finally, it remains to prove the MC when $\beta_{2,t+2}>0$ and the r.h.s. of formula (\ref{diff}) is equal to 0. In this case, since both $\beta_{1,t+1}$ and $\beta_{2,t+1}$ could be 0, we need to show that
\begin{equation}\label{iiii}{t(t+2)(q+3)\over 6}\leq e\leq {t(t+2)(c+3)\over 6}.\end{equation}\indent
But, since we have also proven (\ref{ii}) and (\ref{iii}) when the r.h.s. of (\ref{diff}) is equal to 0, the two inequalities of (\ref{iiii}) have already been shown (respectively, in proving of the first inequality of (\ref{ii}) and the second of (\ref{iii})). This completes the proof of the theorem.{\ }{\ }\qed \\
\\\indent
The above theorem easily generalizes in the following way, giving us the main result of this section:\\
\\\indent
{\bf Theorem 2.7.} {\it Let $A^{'}$ be any codimension 3 algebra whose $h$-vector is the $h$-vector of a compressed algebra. Then the MC holds for $A^{'}$.}\\
\\\indent
{\bf Proof.} Let $A$ be a compressed algebra whose $h$-vector $h$ coincides with that of $A^{'}$. Notice that, given $h$, the socle-vector of $A$ can be uniquely determined - one immediately computes it by induction using the definition of the numbers $r_d$; see Definition-Remark 2.1. Moreover, by the very definition of compressed algebra, we can easily see that the socle-vector of $A^{'}$ must be (entry by entry) greater than or equal to that of $A$.\\\indent
The computations made in the proof of Theorem 2.6 regarding the first two modules of the MFR of $A$ were merely numerical calculations on $h$ (and on the socle-vector of $A$, which is, as we have just said, determined by $h$); from those computations, we found out which shifts had necessarily to appear in the first two modules of the MFR of $A$, and then we used them to prove the bounds of the MC for $A$. Furthermore, since the socle-vector of $A^{'}$ is greater than or equal to that of $A$, in last module of the MFR of $A^{'}$ we must still have the shifts of degree $q+3$ and $c+3$ that we have in that of $A$.\\\indent
Thus, since $A^{'}$ has the $h$-vector of $A$, we can immediately see that the same numerical values we considered in the proof of Theorem 2.6 in bounding the multiplicity of $A$ can be considered when it comes to the algebra $A^{'}$, and therefore we have that the proof of the previous theorem extends to $A^{'}$.{\ }{\ }\qed \\
\\\indent
{\bf Example 2.8.} Consider a codimension 3 algebra $A^{'}=R/I^{'}$, whose associated inverse system module $M^{'}=(I^{'})^{-1}\subset S$ is generated by $L_1^8$, $L_2^8$, $L_3^7$, ..., $L_6^7$, $L_7^6$, ..., $L_{12}^6$, $L_{13}^5$, ..., $L_{20}^5$, where the $L_i$'s are generic linear forms. It is easy to see that the $h$-vector of the algebra $A^{'}$ (which has socle-vector $(0,0,0,0,0,8,6,4,2)$) is $h=(1,3,6,10,15,20,12,6,2)$ (see $[Ia]$, Theorem 4.8 B).\\\indent
A simple calculation shows that $h$ is also the $h$-vector of a compressed (level) algebra (having socle-vector $(0,0,0,0,0,0,0,0,2)$). Therefore, it follows from Theorem 2.7, without any further computations, that $A^{'}$ satisfies the Multiplicity Conjecture.\\
\\
\section{New conjectural bounds for codimension 3 level algebras}
\indent

As we said in the introduction, in general the bounds of the Multiplicity Conjecture seem too loose, although extremely difficult to prove; indeed, it was recently conjectured (see $[HZ]$ and $[MNR2]$) that the only time the lower bound or the upper bound are sharp is when the algebra has a pure MFR. Therefore, we wonder if, in particular cases, we can find {\it ad hoc} sharper bounds, which, possibly, are also easier to handle.\\\indent
In this section we will restrict our attention to codimension 3 level algebras $A$, and construct better candidates for bounding the multiplicity of $A$, which - if true - are in general sharp, even when the MFR of $A$ is not pure. Our bounds have the potentially remarkable advantage to be expressed uniquely in terms of the $h$-vector of $A$, and therefore no knowledge of the MFR of $A$ is required to prove them. Unfortunately, however, the current state of research on codimension 3 level $h$-vectors (e.g., see $[GHMS]$ for a comprehensive overview up to the year 2003, but also our recent surprising results of $[Za4]$) does not yet provide us with a complete picture of which are actually these $h$-vectors and which are not, and therefore we cannot prove our conjectural bounds in general at this point.\\\indent
We will be able to show them, however, in a few interesting particular cases: namely, when the codimension 3 level algebra $A$ is compressed, or when its $h$-vector $h(A)$ ends with $(...,3,2)$. Also, we will prove our lower bound when $h(A)$ begins with $(1,3,h_2,...)$, where $h_2\leq 4$, and our upper bound when $h(A)$ ends with $(...,h_{c-1},h_c)$, where $h_{c-1}\leq h_c+1$.\\\indent
Our bounds, in general, cannot be extended to the non-level case (or to codimension $r>3$). Indeed, our upper bound does not even imply that of the MC as soon as we drop the codimension 3 level hypothesis. Instead, even if our lower bound does not hold for all algebras, we will see that when it is verified for any algebra $A$ of codimension 3, then the lower bound of the MC also holds for $A$. We will exploit this fact to show the lower bound of the MC for any algebra whose $h$-vector begins with $(1,3,4,5,...)$.\\
\\\indent
Let $A=R/I$ be a codimension 3 level algebra, having $h$-vector $h=(1,3,h_2,...,h_c)$ and graded Betti numbers $\beta_{i,j}$, for $i=1,2,3$. The fact that $A$ is level of socle degree $c$ means the the only non-zero number $\beta_{3,j}$ is $\beta_{3,c+3}=h_c$. Also, recall formula (\ref{sh}) (for the codimension 3 case), which provides a relationship between $h$ and the $\beta_{i,j}$'s:
\begin{equation}\label{sh2}h(z)(1-z)^3=1+\sum_{i,j}(-1)^i\beta_{i,j}z^j.\end{equation}\indent
It is easy to see, by (\ref{sh2}), that for any integer $n$, we have
$$\beta_{2,n}-\beta_{3,n}-\beta_{1,n}=h_n-3h_{n-1}+3h_{n-2}-h_{n-3}.$$\indent
Set $$f(n)=h_n-3h_{n-1}+3h_{n-2}-h_{n-3}.$$\indent
Thus, when $f(n)>0$, the Betti number $\beta_{2,n}$ has to be positive, i.e. a shift of degree $n$ must appear in the second module of the MFR of $A$. Conversely, since $A$ is level, if $f(n)<0$ for some $n<c+3$ (when $n=c+3$ we would only be considering $\beta_{3,c+3}$, which is positive since it indicates the socle), then $\beta_{1,n}>0$, i.e. there is a shift of degree $n$ in the first module of the MFR of $A$.\\\indent
Define $i$ as the smallest positive integer such that $f(i)>0$, $j$ as the largest integer such that $f(j)>0$, and $m$ as the largest integer lower than $c+3$ such that $f(m)<0$. As usual, let $t$ be the initial degree of $I$ (that is, the smallest integer such that $f(t)<0$). Then we are ready to state our conjecture:\\
\\\indent
{\bf Conjecture 3.1.} {\it Let $A$ be a codimension 3 level algebra of socle degree $c$ and multiplicity $e$, and let $i$, $j$, $t$ and $m$ be as above. Then
\begin{equation}\label{ho}{t\cdot i\cdot (c+3)\over 6}\leq e\leq {m\cdot j\cdot (c+3)\over 6}.\end{equation}}
\\\indent
From what we observed above in constructing the invariants $i$, $j$ and $m$, we immediately have that our Conjecture 3.1 implies the Multiplicity Conjecture. More precisely:\\
\\\indent
{\bf Proposition 3.2.} {\it Let $A$ be a codimension 3 level algebra satisfying the lower (respectively, upper) bound of Conjecture 3.1. Then $A$ also satisfies the lower (respectively, upper) bound of the MC.}\\
\\\indent
Notice that, with Conjecture 3.1, we are believing in much more than the Multiplicity Conjecture, since the invariants we have constructed to define the bounds of (\ref{ho}) do not always take into account the degrees of the shifts appearing in both of the first two modules of the MFR of $A$ - and these can happen to be many. However, all the examples we know and the computations we have performed so far seem to suggest that the bounds of (\ref{ho}) might hold for all level algebras of codimension 3.\\\indent
In the next example, we will show that the bounds of Conjecture 3.1 in general cannot be improved, even when the MFR of the codimension 3 level algebra $A$ is not pure.\\
\\\indent
{\bf Example 3.3.} Let $A$ be the Gorenstein algebra associated to the inverse system cyclic module generated by $F=y_1^5-y_1y_3^4-y_2^2y_3^3\in S=k[y_1,y_2,y_3]$. It is easy to see (e.g., using $[CoCoA]$, and in particular a program on inverse systems written for us by our friend and colleague Alberto Damiano) that the $h$-vector of $A$ is $h=(1,3,4,4,3,1)$. Hence, we can immediately check that $e=16$, $t=m=2$, $i=j=6$ and $c+3=8$, and therefore that the two bounds of Conjecture 3.1 are sharp for this algebra $A$.\\\indent
Instead, since, as one can easily compute, in the first two modules of the MFR of $A$ there are also shifts of degree 4 (that one cannot notice just by looking at the $h$-vector, since $h_4-3h_3+3h_2-h_1=3-12+12-3=0$), the bounds of the MC are not sharp: indeed, we have $2\cdot 4\cdot 8/6=10.66...<16<4\cdot 6\cdot 8/6=32$.\\
\\\indent
Note that, as soon as we drop the hypothesis that $A$ be level of codimension 3, we are no longer able to recover, from the three invariants $i$, $j$ and $m$ defined above, the same information on the shifts of the MFR of $A$ (or therefore on the MC). In fact, if $r>3$, from (\ref{sh2}) we obtain, in place of $f$, a formula with at least two positive and two negative terms, and therefore the sign of that formula no longer implies the existence of a particular shift in the MFR of $A$.\\\indent
When $r=3$ but we drop the level hypothesis, the same argument holds for the integer $m$. Indeed, the inequality $f(n)<0$ forces $\beta_{3,n}+\beta_{1,n}>0$, but clearly this guarantees neither the positivity of $\beta_{3,n}$ nor that of $\beta_{1,n}$.\\\indent
We just remark here that Francisco ($[Fr]$), who already studied some cases of the Multiplicity Conjecture by looking at possible numerical cancelations among the Betti numbers, suggested an approach (from which he obtained some interesting results) to perform cancelations also when - like in the above cases - there could be more than one way to make them (basically, his choice was to give priority to the rightmost cancelation). However, Francisco's technique and results (which go in a different direction) will not be employed nor further discussed in this paper.\\\indent
Instead, for any codimension 3 algebra $A$, we can see by the same reasoning as above that the invariant $i$ always implies the existence of a shift in the second module of the MFR of $A$. Hence we have:\\
\\\indent
{\bf Proposition 3.4.} {\it Let $A$ be any codimension 3 algebra satisfying the lower bound of (\ref{ho}). Then $A$ also satisfies the lower bound of the MC.}\\
\\\indent
However, the lower bound we supplied in Conjecture 3.1 does not hold for all codimension 3 algebras, as the following example shows:\\
\\\indent
{\bf Example 3.5.} Let $A$ be a codimension 3 compressed algebra having socle-vector $(0,0,0,0,0,0,0,1,0,1)$. Then the $h$-vector of $A$ is $$h=(1,3,6,10,15,21,13,7,3,1).$$\indent
We have $e=80$, $t=6$, $i=7$ and $c+3=12$; therefore, the lower bound of (\ref{ho}) does not hold for $A$, since $${t\cdot i\cdot (c+3)\over 6}=84>80=e.$$\indent
Notice that $A$ is actually an {\it extremely} compressed algebra. In general, for all codimension 3 extremely compressed non-level algebras, it can be shown that neither bound of the MC is sharp (just by relaxing the inequalities \lq \lq $\leq $" into strict inequalities \lq \lq $<$" in our proof of Theorem 2.6 when $q<c$; this result is also consistent with the above-mentioned improvement of the multiplicity conjecture due to Herzog-Zheng and Migliore-Nagel-R\"omer). Thus, since for these algebras the product $t\cdot i\cdot (c+3)/6$ coincides with the upper bound of the MC and is therefore larger that $e$, we have that for all codimension 3 extremely compressed non-level algebras the lower bound of Conjecture 3.1 does not hold.\\
\\\indent
There are several cases, however, when the lower bound of Conjecture 3.1 holds for a non-level algebra $A$. For instance, we can show very easily:\\
\\\indent
{\bf Proposition 3.6.} {\it Let $A$ be any codimension 3 algebra having an $h$-vector which begins with $(1,3,4,5,...)$. Then the lower bound of the MC holds for $A$.}\\
\\\indent
{\bf Proof.} Let $h(A)=(1,h_1,...,h_c)$. By Proposition 3.4 it suffices to show that the first inequality of (\ref{ho}) holds for the algebra $A$ - namely, with the above definitions, that $t\cdot i\cdot (c+3)/ 6\leq e$.\\\indent
It is immediate to check that $t=2$ and $i=3$. Hence what we want to prove is that $e\geq c+3$. In fact, for any $c\geq 3$, we clearly have that $e\geq 1+3+4+5+(c-3)=c+10$, and the result follows.{\ }{\ }\qed \\
\\\indent
Let us now come back to the codimension 3 level case. As we said, we are not yet able to prove Conjecture 3.1, and {\it a fortiori} the Multiplicity Conjecture, in full generality, but there are a few interesting cases where we can be successful. Namely, we have:\\
\\\indent
{\bf Theorem 3.7.}. {\it Let $A$ be a codimension 3 level algebra, and let $h=(1,3,h_2,...,h_{c-1},h_c)$ be its $h$-vector. Then:\\\indent
In the following cases Conjecture 3.1 holds for $A$:\\\indent
i). $A$ is compressed.\\\indent
ii). $h_{c-1}=3$ and $h_c=2$.\\\indent
In the following case the lower bound of Conjecture 3.1 holds for $A$:\\\indent
iii). $h_2\leq 4$.\\\indent
In the following case the upper bound of Conjecture 3.1 holds for $A$:\\\indent
iv). $h_{c-1}\leq h_c+1$.}\\
\\\indent
{\bf Proof.} i). This case is already implicitly shown by the argument of Theorem 2.6 (if we let $q=c$, i.e. we consider that $A$ be level): indeed, the bounds of Conjecture 3.1 are exactly those we have shown in that proof, since the shifts we have considered in the first two modules of the MFR of $A$ in order to prove the MC were always those whose existence was forced by inequalities on the Betti numbers, and we have used exactly the same inequalities to construct the invariants $i$, $j$ and $m$. This proves the theorem when $A$ is compressed.\\\indent
ii). In $[Za3]$, Theorem 2.9, we characterized the level $h$-vectors of the form $(1,3,...,3,2)$ as those which can be expressed as the sum of  $(0,1,1,...,1)$ and a codimension 2 Gorenstein $h$-vector. Since the latter $h$-vectors are of the form $(1,2,...,p-2,p-1,p-1,...,p-1,p-2,...,2,1)$ for some integer $p\geq 3$ (this fact is easy to see, and was first noticed by Macaulay in $[Ma]$), we have that the codimension 3 level $h$-vectors we are considering here are of the form either $h=(1,3,3,...,3,3,2)$ or $h=(1,3,4,...,p,p,...,p,p-1,...,4,3,2)$.\\\indent
Hence, if $c$ as usual denotes the socle degree of $h(A)$, $c\geq 2$, it is easy to compute that $h$ stabilizes at $p$ exactly $c-2(p-2)+1=c-2p+5$ times. In particular, $3\leq p\leq (c+4)/2$.\\\indent
The multiplicity of $A$ is $$e=(1+3+4+...+p-1)+p(c-2p+5)+(p-1+p-2+...+3+2)=$$$${(p-1)p\over 2}-2+p(c-2p+5)+{(p-1)p\over 2}-1=pc-2p^2+5p+p^2-p-3=pc-(p-3)(p-1).$$\indent
Let $p=3$, i.e. $h=(1,3,3,...,3,2)$. Then straightforward computations show that $t=2$, $i=3$, $j=c+2$, $m=c$ for $c\neq 3$ and $m=2$ for $c=3$. Checking the bounds of Conjecture 3.1, i.e. proving that, under the current conditions, $${t\cdot i\cdot (c+3)\over 6}\leq pc-(p-3)(p-1)\leq {m\cdot j\cdot (c+3)\over 6},$$ is an easy exercise which is left to the reader.\\\indent
Let $p=4$. We have $t=2$, $i=4$ for $c\neq 5$ and $i=5$ for $c=5$, $j=c+2$, $m=c+1$. Again, it is trivial to check the bounds of Conjecture 3.1 (just notice that here $c\geq 4$, since $p=4$).\\\indent
Finally, let $p\geq 5$. Thus, $t=2$, $i=3$, $j=c+2$, $m=c+1$. We want to show that \begin{equation}\label{pp}{2\cdot 3\cdot (c+3)\over 6}\leq pc-(p-3)(p-1)\leq {(c+1)(c+2)(c+3)\over 6}.\end{equation}\indent
The first inequality of (\ref{pp}) is equivalent to $c+3\leq pc-(p-3)(p-1)$, i.e. to $(p-1)(c-p+3)\geq 3$, and this is true since $p-1\geq 4$ and $c-p+3\geq (2p-4)-p+3=p-1\geq 4$.\\\indent
As for the second inequality of (\ref{pp}), it is clearly enough to show that $pc\leq (c+1)(c+2)(c+3)/6$, or therefore that $p\leq (c+2)(c+3)/6$. But the latter inequality is easily verified, since $p\leq (c+4)/2$ and $c\geq 6$ for $p\geq 5$.\\\indent
This completes the proof of (\ref{pp}), and that of case ii) of the theorem.\\\indent
iii). Using the tables for the possible codimension 3 level $h$-vectors of low socle degree (see $[GHMS]$, Appendix F), and employing a few standard considerations, one can determine the form of the level $h$-vectors beginning with $(1,3,h_2,...)$, for $h_2\leq 4$. Precisely:\\\indent
For $h_2=1$, $h$ can only be $(1,3,1)$.\\\indent
For $h_2=2$, the only possibility for $h$ is $(1,3,2)$ (since $h_3$ must be lower than 3 by Macaulay's theorem, cannot be 2 otherwise would force a one-dimensional socle in degree 1 (e.g., see $[Za5]$, Theorem 3.5), and cannot be 1 because of the symmetry of Gorenstein $h$-vectors. Thus, $h_3=0$).\\\indent
For $h_2=3$, by considerations similar to those we have made for the previous case, we see that $h$ must necessarily be of the form $(1,3,3,...,3,h_c)$, where $h_c$ equals 1, 2 or 3.\\\indent
Let $h_2=4$. Then, likewise, we can show that all the possibilities for $h$ are: $(1,3,4,5,...)$, $(1,3,4,4,4,u,...)$ (where $u\leq 4$), $(1,3,4,4,3,2)$, $(1,3,4,4,3,1)$, $(1,3,4,4,3)$, $(1,3,4,4,2)$, $(1,3,4,4)$, $(1,3,4,3,2)$, $(1,3,4,3,1)$, $(1,3,4,3)$, $(1,3,4,2)$, $(1,3,4)$.\\\indent
In all of the cases listed above it is easy to check that the lower bound of Conjecture 3.1 is verified, so we will avoid the computations here. This completes the proof of this case.\\\indent
iv). From the hypothesis $h_{c-1}\leq h_c+1$, it clearly follows that $-h_{c-1}>-3h_c$, and therefore $j=c+2$. First suppose that $h_{c-1}\leq h_c$. Then, since $-h_{c-2}<0\leq 3(h_c-h_{c-1})$, we have $m=c+1$. Hence the upper bound of Conjecture 3.1 is $(c+1)(c+2)(c+3)/6$. We have $$e\leq 1+3+...+N(3,c)=N(4,c)=(c+1)(c+2)(c+3)/6,$$
as we wanted to show.\\\indent
Now let $h_{c-1}=h_c+1$. It is immediate to check, by definition of $m$, that this time we have $m=c+1$ if and only if $h_{c-2}>3$. So, when $h_{c-2}>3$, we are done by reasoning as above. Instead, by using the same standard techniques we have employed in showing point iii) of this theorem, we can see that all the codimension 3 level $h$-vectors such that $h_{c-1}=h_c+1$ and $h_{c-2}\leq 3$ are: $(1,3,3,...,3,2)$ (for any $c\geq 2$), $(1,3,4,3)$, $(1,3,5,4)$ and $(1,3,6,5)$. It is easy to check the upper bound of Conjecture 3.1 for these $h$-vectors (for the first $h$-vector see the case $p=3$ of point ii) above), so we will leave the computations as an exercise for the reader. This completes the proof of this case and that of the theorem.{\ }{\ }\qed \\
\\\indent
{\bf Acknowledgements.} We wish to thank the referee for her/his insightful comments, which helped improve the exposition of this paper.\\\indent
This author is funded by the G\"oran Gustafsson Foundation.\\
\\
\\
\\
{\bf \huge References}\\
\\
$[CoCoA]$ A. Capani, G. Niesi and L. Robbiano: {\it CoCoA, a system for doing computations in commutative algebra}, available via anonymous ftp, cocoa.dima.unige.it.\\
$[EI]$ {\ } J. Emsalem and A. Iarrobino: {\it Some zero-dimensional generic singularities; finite algebras having small tangent space}, Compositio Math. 36 (1978), 145-188.\\
$[Fr]$ {\ } C. Francisco: {\it New approaches to bounding the multiplicity of an ideal}, J. of Algebra, to appear (preprint: math.AC/0506024).\\
$[FS]$ {\ } C. Francisco and H. Srinivasan: {\it Multiplicity conjectures}, preprint (2005).\\
$[FL]$ {\ } R. Fr\"oberg and D. Laksov: {\it Compressed Algebras}, Conference on Complete Intersections in Acireale, Lecture Notes in Mathematics, No. 1092 (1984), 121-151, Springer-Verlag.\\
$[Ge]$ {\ } A.V. Geramita: {\it Inverse Systems of Fat Points: Waring's Problem, Secant Varieties and Veronese Varieties and Parametric Spaces of Gorenstein Ideals}, Queen's Papers in Pure and Applied Mathematics, No. 102, The Curves Seminar at Queen's (1996), Vol. X, 3-114.\\
$[GHMS]$ {\ } A.V. Geramita, T. Harima, J. Migliore and Y.S. Shin: {\it The Hilbert Function of a Level Algebra}, Memoirs of the Amer. Math. Soc., to appear.\\
$[GV]$ {\ } E. Guardo and A. Van Tuyl: {\it Powers of complete intersections: graded Betti numbers and applications}, Illinois J. of Math. 49 (2005), 265-279.\\
$[HS]$ {\ } J. Herzog and H. Srinivasan: {\it Bounds for multiplicities}, Trans. Amer. Math. Soc. 350 (1998), 2879-2902.\\
$[HZ]$ {\ } J. Herzog and X. Zheng: {\it Notes on the multiplicity conjecture}, Collectanea Math. 57 (2006), 211-226.\\
$[HM]$ {\ } C. Huneke and M. Miller: {\it A note on the multiplicity of Cohen-Macaulay algebras with pure resolutions}, Canad. J. Math 37 (1985), 1149-1162.\\
$[Ia]$ {\ } A. Iarrobino: {\it Compressed Algebras: Artin algebras having given socle degrees and maximal length}, Trans. Amer. Math. Soc. 285 (1984), 337-378.\\
$[IK]$ {\ } A. Iarrobino and V. Kanev: {\it Power Sums, Gorenstein Algebras, and Determinantal Loci}, Springer Lecture Notes in Mathematics (1999), No. 1721, Springer, Heidelberg.\\
$[Ma]$ {\ } F.H.S. Macaulay: {\it The Algebraic Theory of Modular Systems}, Cambridge Univ. Press, Cambridge, U.K. (1916).\\
$[MNR1]$ {\ } J. Migliore, U. Nagel and T. R\"omer: {\it The Multiplicity Conjecture in low codimension}, Math. Res. Letters 12 (2005), 731-748.\\
$[MNR2]$ {\ } J. Migliore, U. Nagel and T. R\"omer: {\it Extensions of the Multiplicity Conjecture}, Trans. Amer. Math. Soc., to appear (preprint: math.AC/0505229).\\
$[Mi]$ {\ } R. Mir\'o-Roig: {\it A note on the multiplicity of determinantal ideals}, J. of Algebra, to appear (preprint: math.AC/0504077).\\
$[Ro]$ {\ } T. R\"omer: {\it Note on bounds for multiplicities}, J. of Pure and Applied Algebra 195 (2005), 113-123.\\
$[Za1]$ {\ } F. Zanello: {\it Extending the idea of compressed algebra to arbitrary socle-vectors}, J. of Algebra 270 (2003), No. 1, 181-198.\\
$[Za2]$ {\ } F. Zanello: {\it Extending the idea of compressed algebra to arbitrary socle-vectors, II: cases of non-existence}, J. of Algebra 275 (2004), No. 2, 730-748.\\
$[Za3]$ {\ } F. Zanello: {\it Level algebras of type 2}, Comm. in Algebra 34 (2006), No. 2, 691-714.\\
$[Za4]$ {\ } F. Zanello: {\it A non-unimodal codimension 3 level $h$-vector}, preprint (math.AC/0505678).\\
$[Za5]$ {\ } F. Zanello: {\it When is there a unique socle-vector associated to a given $h$-vector?}, Comm. in Algebra, to appear (preprint: math.AC/0411229).

}

\end{document}